\documentclass[a4paper,11pt,reqno]{amsart}
\usepackage{amsmath,a4,amsfonts,amssymb,amsthm,enumerate,fleqn,breqn,dsfont}
\usepackage{enumitem}
\setdescription{leftmargin=0cm,labelindent=0cm}
\usepackage[a4paper,left=1in, right=1in, top=1.4in, bottom=1.1in]{geometry}

\usepackage{etoolbox}
\patchcmd{\section}{\scshape}{\bfseries}{}{}
\makeatletter
\renewcommand{\@secnumfont}{\bfseries}
\makeatother

\def\RR{\mathbb{R}}
\def\NN{\mathbb{N}}
\def\HH{\mathfrak{Q}}
\def\U{\{u_i\}_{i\in I}}

\def\T{\mathfrak{T}}
\def\cc#1{\{#1\}}
\def\pp#1{\|#1\|}
\def\ra{\rangle}
\def\la{\langle}
\def\nh{\mathcal{H}}

\def\E{\mathcal{E}}
\def\H{V_R(\HH)}

\def\E{\ell_2(\HH)}
\def\cc#1{\{#1\}}
\def\pp#1{\|#1\|}
\def\ra{\rangle}
\def\la{\langle}

\theoremstyle{plain}
\newtheorem{theorem}{\bf Theorem}[section]

\theoremstyle{remark}
\newtheorem{definition}[theorem]{\bf Definition}

\newtheorem{example}[theorem]{\bf Example}

\title[Duals of a frame in quaternionic Hilbert Spaces]{Duals of a frame in quaternionic Hilbert Spaces}
\author[ S\lowercase{harma}, S\lowercase{ingh} \lowercase{and} S\lowercase{ahu}]
{S.K. S\lowercase{harma}$^\dag$, G\lowercase{hanshyam} S\lowercase{ingh}$^\S$ \lowercase{and} S\lowercase{oniya} S\lowercase{ahu}$^\ddag$ \bigskip \\
	$^{\dag}$D\lowercase{epartment} \lowercase{of} M\lowercase{athematics},\\
	K\lowercase{irori} M\lowercase{al} C\lowercase{ollege},\\
	U\lowercase{niversity of} D\lowercase{elhi}, D\lowercase{elhi~110~007}, INDIA.\\
	\emph{$^\dag$E}\lowercase{\emph{-mail}: sumitkumarsharma@gmail.com}\bigskip\\
	$^{\ddag, \S}$D\lowercase{epartment} \lowercase{of} M\lowercase{athematics and}  S\lowercase{tatistics},\\
	M.L.S U\lowercase{niversity}, U\lowercase{daipur}, INDIA.\\
	\emph{$^\S$E}\lowercase{\emph{-mail}: ghanshyamsrathore@yahoo.co.in}\\
\emph{$^\ddag$E}\lowercase{\emph{-mail}: soniya.vvsjay1@gmail.com}}


\begin{document}
	\maketitle \baselineskip14pt



 \baselineskip12pt
\begin{abstract}
 Frames in a separable quaternionic Hilbert space were introduced and studied in \cite{SS}   to have more applications. In this paper, we extend the study of frames in quaternionic Hilbert
  spaces and introduce  different types of duals of a frame in separable quaternionic Hilbert spaces. As   an application, we give the orthogonal projection of $\ell^2(\HH)$ onto the range of analysis operator of the given frame, in terms of elements of  canonical dual frame and  elements of the  frame in quaternionic Hilbert space.  Finally, we  give an expression for the orthogonal projection in terms of operators related to the frame and its canonical dual frame in quaternionic Hilbert space.
 \end{abstract}\baselineskip14pt

\section{Introduction}
\def\xmn{\cc{u_i}_{i\in \NN}}
\def\hmn{\cc{\nh_{n,i}}_{i=1,2, \cdots, m_n \atop{n\in \NN}}}
\def\h2n{\cc{\nh_{n,i}}_{i=1,2, \cdots, 2n \atop{n\in \NN}}}
\def\tmn{\cc{\mathfrak{T}_{n,i}}_{i=1,2, \cdots, m_n \atop{n\in \NN}}}
\def\Rmn{\cc{\mathfrak{R}_{n,i}: \nh \to \nh_{n,i}}_{i=1,2, \cdots, m_n \atop{n\in \NN}}}
\def\Tmn{\cc{\mathfrak{T}_{n,i}: \nh \to \nh_{n,i}}_{i=1,2, \cdots, m_n \atop{n\in \NN}}}
\def\Umn{\cc{\mathfrak{U}_{n,i}: \nh \to \nh_{n,i}}_{i=1,2, \cdots, m_n \atop{n\in \NN}}}
\def\Tnn{\cc{\mathfrak{T}_{n,i}: \nh \to \nh_{n,i}}_{i=1,2, \cdots, 2n \atop{n\in \NN}}}
\def\ymn{\cc{v_i}_{i\in\NN}}
\def\xxmn{\cc{\la x, x_{n,i}\ra}_{i=1,2, \cdots, m_n \atop{n\in \NN}}}
\def\alphamn{\cc{\alpha_{n,i}}_{i=1,2, \cdots, m_n \atop{n\in \NN}}}
\def\betamn{\cc{\beta_{n,i}}_{i=1,2, \cdots, m_n \atop{n\in \NN}}}
\def\llim{\lim\limits_{n\to \infty}}
\def\S{{S}}
\def\T{\mathfrak{T}_{\cc{m_n}}}
\def\suml{\sum\limits_{i=1}^{\infty}}
\def\mnsum{\sum\limits_{i=1}^{\infty}}
\def\mpsum{\sum\limits_{i=1}^{m_p}}
\def\mqsum{\sum\limits_{i=1}^{m_q}}
\def\elii{\ell^2(\HH)}

Formally, frames for Hilbert spaces (in particular for $L^2[a,b]$) were introduced way back in  1952 by Duffin and Schaeffer \cite{DS}
 as a tool to study of non-harmonic Fourier series. They defined the following

``A sequence $\{x_n\}_{n\in\NN}$ in a Hilbert space $\nh$ is said to be a \emph{frame}
for $\nh$ if there exist constants $A$ and $B$ with $0<A\le B<\infty$
such that
\begin{eqnarray}
A\|x\|^2\le \sum\limits_{n=1}^\infty |\langle
x,x_n\rangle|^2\le B\|x\|^2, \ \ \text{for all} \ x\in \nh."
\end{eqnarray}
Moreover, the positive constants $A$ and $B$, respectively, are called \textit{lower}
and\textit{ upper} frame bounds for the frame $\cc{x_n}_{n\in\NN}$. The inequality
$(1.1)$ is called the \emph{frame inequality} for the frame
$\{x_n\}_{n\in\NN}$. A sequence $\cc{x_n}_{n\in\NN} \subset \nh$ is called a \emph{Bessel sequence} if it satisfies
upper frame inequality in $(1.1)$.
A frame $\cc{x_n}_{n\in\NN}$ in $\nh$ is said to be
\begin{itemize}[leftmargin=.5in]
\item \emph{{tight}} if it is possible to choose  $A,\ B$ satisfying inequality $(1.1)$ with $A=B$.
\item \emph{Parseval} if it is possible to choose  $A,\ B$ satisfying inequality $(1.1)$ with $A=B=1$.
\item \emph{exact} if removal of any  $x_n$ renders the collection $\cc{x_i}_{i\ne n}$ no longer a frame for $\nh$.
\end{itemize}

\medskip

 Later, frames were further reintroduced, in 1986 by, Daubechies, Grossmann
and Meyer \cite{DGM}, they observed that frames can be used to approximate functions in $L^2({\RR})$. One can also considered
frames as one of the generalizations of orthonormal bases in Hilbert spaces and being  redundant
  frames expansions are more useful and advantageous over basis expansions in a variety of practical
   applications. Now a days, frames are regarded as
one of an important tool  to study various areas like representation
of signals, characterization of function spaces and other fields
of applications such as: signal and image processing \cite{CD}, filter bank theory \cite{BHF},
wireless communications \cite{HP} and sigma-delta quantization \cite{BPY}. For more literature on
frame theory, one may refer to \cite{C1,CH2}.

In recent years, many generalizations of frames have been introduced and studied.
In 2004, Casazza and Kutyniok \cite{CK} defined  frames of subspaces
(frames of subspaces has many applications  in sensor networks and packet encoding),
 Li and Ogawa \cite{LO} introduced the notion of pseudo-frames in Hilbert spaces using
 two Bessel sequences, Fornasier \cite{F} introduced the notion of bounded quasi-projectors,
 Christensen and Eldar \cite{CE} gave  oblique frames. In 2006, Sun \cite{SUN} introduced  generalized frames or
$g$-frames for Hilbert spaces and proved that frames of subspaces (fusion frames), pseudo frames,
bounded quasi-projectors and oblique frames are special cases of $g-$frames.

Recently, Khokulan, Thirulogasanthar and Srisatkunarajah \cite{KTS} introduced and studied frames for finite dimensional quaternionic Hilbert spaces. Sharma and Virender \cite{SV} study some different types of dual frames of a given frame in a finite dimensional   quaternionic Hilbert space and gave various types of reconstructions with the help of  dual frame. Very recently, Sharma and Goel \cite{SS} introduced and studied frames for seperable quaternionic Hilbert spaces  and  Chen, Dang and Qian \cite{CDQ} had studied frames for Hardy
spaces in the contexts of the quaternionic space and the Euclidean space in the Clifford
algebra.

In this paper,  we extend the study of frames in quaternionic Hilbert
spaces and introduce  different types of duals of  a frame in separable quaternionic Hilbert spaces. As   an application, we give the orthogonal projection of $\ell^2(\HH)$ onto the range of analysis operator of the given frame, in terms of elements of  canonical dual frame and elements of the  frame  in quaternionic Hilbert  spaces.  Finally, we  give an expression for the orthogonal projection in terms of operators related to the frame and its canonical dual frame in quaternionic Hilbert spaces.

\def\nh{V_R(\HH)}
	
\section{Quaternionic Hilbert space}
\setcounter{equation}{0}
\fontsize{12}{14}

As the quaternions are non-commutative in nature therefore there are two different types of quaternionic Hilbert spaces, the left  quaternionic Hilbert space and the right quaternionic Hilbert space depending on positions of quaternions. In this section, we will study some basic notations about the algebra of quaternions, right quaternionic Hilbert spaces and operators on right quaternionic Hilbert spaces. 

Throughout this paper, we will denote $\mathfrak{Q}$ to be a non-commutative field of quaternions,  $I$ be a non empty countable set of indicies, $\H$ be a separable right quaternionic Hilbert space,  by the	term ``right linear operator", we mean a ``right $\HH$-linear operator" and $\mathfrak{B}(\H)$ denotes the set of all bounded (right $\HH$-linear) operators of $\H$:
\begin{eqnarray*}
\mathfrak{B}(\H) := \{T : \H\rightarrow\H: \|T\|<\infty\}.
\end{eqnarray*}

The non-commutative field of quaternions $\mathfrak{Q}$ is a four dimensional real algebra with unity. In $\HH$, $0$ denotes the null element and $1$ denotes the identity with respect to multiplication. It also includes three so-called imaginary units, denoted by $i,j,k$. i.e.,
\begin{eqnarray*}
\mathfrak{Q}=\cc{x_0+x_1i +x_2j +x_3k \ :\ x_0,\ x_1,\ x_2,\  x_3\in \RR}
\end{eqnarray*}
where $i^2=j^2=k^2=-1; \ ij=-ji=k; \ jk=-kj=i$ and $ki=-ik=j$. For each quaternion $q=x_0+x_1i +x_2j +x_3k \in \mathfrak{Q}$, the  conjugate of $q$ is
denoted by $\overline{q}$ and is defined as
\begin{eqnarray*}
\overline{q}=x_0-x_1i -x_2j -x_3k \in \mathfrak{Q}.
\end{eqnarray*}
If $q=x_0+x_1i +x_2j +x_3k$ is a quaternion, then $x_0$ is called the real part of $q$ and $x_1i +x_2j +x_3k$ is called the imaginary part  of $q$. The modulus of a quaternion $q=x_0+x_1i +x_2j +x_3k$ is defined as
\begin{eqnarray*}
|q|=(\overline{q}q)^{1/2} = (q\overline{q})^{1/2}= \sqrt{x_0^2 +x_1^2 +x_2^2 +x_3^2 }.
\end{eqnarray*}
For every non-zero quaternion $q=x_0+x_1i +x_2j +x_3k \in \mathfrak{Q}$, there exists a unique inverse $q^{-1}$ in $\mathfrak{Q}$ as
\begin{eqnarray*}
q^{-1}=\dfrac{\overline{q}}{|q|^2 } = \dfrac{x_0-x_1i -x_2j -x_3k }{{x_0^2 +x_1^2 +x_2^2 +x_3^2 }}.
\end{eqnarray*}

\begin{definition}
	A \textit{right quaternionic vector space} $\mathds{V}_R(\HH)$ is a   vector space under right scalar multiplication over the field of quaternionic $\HH$, i.e.,
	\begin{eqnarray}\label{2.1}
	\mathds{V}_R(\HH)\times\HH &\rightarrow& \mathds{V}_R(\HH) \nonumber\\ 
	(u,q)&\rightarrow& uq
	\end{eqnarray}
	and for each $u, v\in\mathds{V}_R(\HH)$ and $p, q\in\HH$, the right scalar multiplication (\ref{2.1}) satisfying the following properties: 
	\begin{eqnarray*}
	&&(u+v)q=uq+vq\\
	&&u(p+q)=up+uq\\
	&&v(pq)=(vp)q.
	\end{eqnarray*}
\end{definition}

\begin{definition}\label{qhs}
	A \textit{right quaternionic pre-Hilbert space} or \textit{right quaternionic inner product space} $\mathds{V}_R(\HH)$ is a right quaternionic vector space together with the binary mapping
	$\langle . | . \rangle : \mathds{V}_R(\HH) \times \mathds{V}_R(\HH) \to \mathfrak{Q}$ (called the \textit{Hermitian quaternionic inner product})
	which satisfies following properties:
	\begin{enumerate}[label=(\alph*)]
		\item $\overline{\langle v_1 | v_2 \rangle} = \langle v_2 | v_1 \rangle$  for all $v_1, v_2 \in \mathds{V}_R(\HH)$.
		\item $\langle v | v \rangle > 0 \ \text{for all} \ \ 0 \ne  v\in \nh.$
		\item $\langle v | v_1 + v_2 \rangle = \langle v | v_1 \rangle + \langle v | v_2 \rangle $  for all $v, v_1, v_2 \in \mathds{V}_R(\HH).$
		\item $\langle v | uq \rangle = \langle v | u \rangle q $  for all $v, u \in \mathds{V}_R(\HH)$ and $ q \in \mathfrak{Q}$.
	\end{enumerate}
\end{definition}

Let $\mathds{V}_R(\HH)$ be right quaternionic inner product space with the Hermitian inner product $\langle .|.\rangle$. ~Define the quaternionic norm $\|.\|:\mathds{V}_R(\HH)\rightarrow\RR^+$ on $\mathds{V}_R(\HH)$ by
\begin{eqnarray}\label{2.2}
\|u\|=\sqrt{\langle u|u\rangle},\ u\in\mathds{V}_R(\HH).
\end{eqnarray}

\begin{definition}
	The right quaternionic pre-Hilbert space  is called a \textit{right quaternionic Hilbert space}, if it is complete with respect to the norm (2.2) and is denoted by $\H$.
\end{definition}

\begin{theorem}[The Cauchy-Schwarz Inequality]\cite{GMP} 
	If $\H$ is a right quaternionic Hilbert space then  
	\begin{eqnarray*}
	|\langle u|v\rangle|^2\leq\langle u|u\rangle\langle v|v\rangle, \ \ \text{for all}\ \ u, v \in \H.
	\end{eqnarray*}
	Moreover, a norm as defined in (\ref{2.2}) satisfies the following properties:
	\begin{enumerate}[label=(\alph*)]
		\item $\|uq\|=\|u\||q|$, for all $u\in\H$ and $q\in\HH$.
		\item $\|u+v\|\leq\|u\|+\|v\|$, for all $u, v\in\H$.
		\item $\|u\|=0$ for some $u\in\H$, then $u=0$.
	\end{enumerate}
\end{theorem}

\noindent
For the non-commutative field of quaternions $\HH$,  define the quaternionic Hilbert space $\ell_2(\HH)$ by
\begin{eqnarray*}
\ell_2(\HH) = \bigg\{\{q_i\}_{i\in \NN}\subset \HH : \ \sum_{i\in \NN}|q_i|^2<+\infty \bigg\}
\end{eqnarray*}
under right multiplication by quaternionic scalars together with   the  quaternionic inner product on $\ell_2(\HH)$ defined as
\begin{eqnarray}\label{2.3}
\langle p|q\rangle=\sum_{i\in \NN}\overline{p_i}q_i,\ p=\{p_i\}_{i\in \NN}\ \text{and}\ q=\{q_i\}_{i\in \NN}\in\ell_2(\HH).
\end{eqnarray}
It is easy to observe that $\ell_2(\HH)$ is a right quaternionic Hilbert space with respect to quaternionic inner product (\ref{2.3}).

\begin{definition}[\cite{GMP}] Let $\H$ be a right quaternionic Hilbert Space and  $S$ be a subset of $\H$. Then, define the set:
	\begin{itemize}[leftmargin=.35in]
		\item	$S^{\bot}=\{v\in\H:\langle v|u\rangle=0\ \forall\ u\in S\}.$
		\item $\langle S\rangle$ be the right $\HH$-linear subspace of $\H$ consisting of all finite right $\HH$-linear combinations of elements of $S$.
	\end{itemize}
\end{definition}

\begin{theorem}[\cite{GMP}]\label{2.6t}
	Let $\H$ be a quaternionic Hilbert space and let $N$ be a subset of $\H$
	such that, for $z, z'\in N$ such that $\langle z|z'\rangle=0$ if $z\neq z'$ and $\langle z|z\rangle=1$. Then the following conditions are equivalent:
	\begin{enumerate}[label=(\alph*)]
		\item For every $u,v\in\H$, the series $\sum_{z\in N}\langle u|z\rangle\langle z|v\rangle$ converges absolutely and
		\begin{eqnarray*}
		\langle u|v\rangle=\sum_{z\in N}\langle u|z\rangle\langle z|v\rangle.
		\end{eqnarray*}
		\item For every $u\in\H$, $\|u\|^2=\sum\limits_{z\in N}|\langle z|u\rangle|^2$.
		\item $N^\bot={0}$.
		\item $\langle N\rangle$ is dense in $H$.
	\end{enumerate}
\end{theorem}

\begin{definition}[\cite{GMP}]
	Every quaternionic Hilbert space $\H$ admits a subset $N$, called \textit{Hilbert basis or orthonormal basis} of $\H$, such that, for $z, z'\in N$, $\langle z|z'\rangle=0$ if $z\neq z'$ and $\langle z|z\rangle=1$ and satisfies all the conditions of Theorem \ref{2.6t}.
\end{definition}
Further, if there are two such sets, then they have the same cardinality. Furthermore, if $N$ is a Hilbert basis of $\H$, then for every $u\in\H$ can be uniquely expressed as
\begin{eqnarray*}
u=\sum_{z\in N}z\langle z|u\rangle
\end{eqnarray*}
where the series $\sum\limits_{z\in N}z\langle z|u\rangle$ converges absolutely in $\H$.

\begin{definition}[\cite{adler}]
	Let $\H$ be a right quaternionic Hilbert space and $T$ be an operator on $\H$. Then $T$ is said to be
	\begin{itemize}[leftmargin=.25in]
		\item \emph{right $\HH$-linear} if
		$T( v_1\alpha + v_2\beta) = T(v_1)\alpha +  T(v_2)\beta, \ \text{for all} \ v_1, v_2 \in \H \ \text{and} \ \alpha, \beta \in \HH.$
		\item \emph{bounded} if there
		exist $K\ge 0$ such that
		$\pp{T(v)} \le K \pp{v}, \ \text{for all} \ v\in \H.$
	\end{itemize}
\end{definition}

\begin{definition}[\cite{adler}]
	Let $\H$ be a right quaternionic Hilbert space and $T$ be an operator on $\H$. Then the \textit{adjoint operator} $T^*$ of $T$ is defined by
	\begin{eqnarray*}
	\la v|Tu\ra = \la T^*v|u\ra, \ \text{for all} \ u, v \in\H.
	\end{eqnarray*}
	Further, $T$ is said to be \emph{self-adjoint} if $T=T^*$.
\end{definition}

\begin{theorem}[\cite{adler}]
	Let $\H$ be a right quaternionic Hilbert space and $S$ and  $T$ be two bounded right $\HH$-4linear operators on $\H$. Then
	\begin{enumerate}[label=(\alph*)]
		\item $T+S$ and $TS\in\mathfrak{B}(\H)$. Moreover:
		\begin{eqnarray*}
		\|T+S\|\leq\|T\|+\|S\|\ \text{and}\ \|TS\|\leq\|T\|\|S\|.
		\end{eqnarray*}
		\item $\la Tv|u\ra = \la v| T^*u \ra$.
		\item $(T+S)^*=T^*+S^*$.
		\item $(TS)^*=S^*T^*$.
		\item $(T^*)^*=T.$
		\item $I^*=I$, where $I$ is the identity operator on $\H$.
		\item If $T$ is an invertible operator then $(T^{-1})^*=(T^*)^{-1}$.
	\end{enumerate}
\end{theorem}


\section{Frames in quaternionic Hilbert spaces}

We begin this section with the following definition of frames in a  separable right quaternionic Hilbert space $\H$ defined in \cite{SS}:
\begin{definition}Let $\H$ be a right quaternionic Hilbert space and $\{u_i\}_{i\in \NN}$ be a sequence in $\H$. Then $\cc{u_i}_{i\in \NN}$ is said to be a \textit{frame} for $V_R(\HH)$, if there exist two finite constants with $0<A\le B$  such that
	\begin{eqnarray}\label{3.1}
	A\|u\|^2\leq\sum_{i\in \NN}|\langle u_i|u\rangle|^2\leq B\|u\|^2, \ \text{for all}\ u\in V_R(\HH).
	\end{eqnarray}
	The positive constants $A$ and $B$, respectively, are called \textit{lower}  and \textit{upper} frame bounds for the frame $\{u_i\}_{i\in \NN}$. The inequality (\ref{3.1}) is called \textit{frame inequality} for the frame $\{u_i\}_{i\in I}$. A sequence $\{u_i\}_{i\in \NN}$ is called a \textit{Bessel sequence} for a right quaternionic Hilbert space $\H$ with bound $B$, if $\{u_i\}_{i\in \NN}$ satisfies the right hand side of the inequality (\ref{3.1}).
	A frame $\{u_i\}_{i\in \NN}$ for a right quaternionic Hilbert space $V_R(\HH)$ is said to be
	\begin{itemize}[leftmargin=.25in]
		\item \textit{tight}, if it is possible to choose $A$ and $B$ satisfying inequality (3.1) with $A=B$.
		\item \textit{Parseval frame}, if it is tight with $A=B=1$.
		\item \textit{exact}, if it ceases to be a frame whenever anyone of its element is removed.
	\end{itemize}
\end{definition}

 If $\{u_i\}_{i\in \NN}$ is a Bessel sequence for a right quaternionic Hilbert space $V_R(\HH)$. Then, the \textit{(right) synthesis operator} for $\{u_i\}_{i\in \NN}$ is a right linear operator $T:\ell_2(\HH)\to V_R(\HH)$ defined by
\begin{eqnarray*}T(\{q_i\}_{i\in \NN})=\sum_{i\in \NN}u_iq_i,\ \ \{q_i\}_{i\in \NN}\in\ell_2(\HH).
\end{eqnarray*}
The adjoint operator $T^*$ of right synthesis operator $T$ is called the \textit{(right) analysis operator}. Further,  the  analysis operator $T^*:V_R(\HH)\to \ell_2(\HH)$  is given by
\begin{eqnarray*}
T^*(u)=\{\langle u_i|u\rangle\}_{i\in \NN},\ u\in V_R(\HH).
\end{eqnarray*}
Infact, for $u\in V_R(\HH)$ and $\{q_i\}_{i\in \NN}\in\ell_2(\HH)$, we have
\begin{eqnarray*}
\langle T^*(u)|\{q_i\}_{i\in \NN}\rangle&=&\langle u|T(\{q_i\}_{i\in \NN})\rangle\\
&=&\bigg\langle u\bigg|\sum_{i\in \NN}u_iq_i\bigg\rangle\\
&=&\sum_{i\in\NN}\langle u|u_i\rangle q_i\\
&=&\bigg\langle \{\langle u_i|u\rangle\}_{i\in \NN}, \{q_i\}_{i\in \NN}\bigg\rangle.
\end{eqnarray*}
Thus
\begin{eqnarray*}
T^*(u)=\{\langle u_i|u\rangle\}_{i\in \NN},\ u\in V_R(\HH).
\end{eqnarray*}

\begin{theorem}[\cite{SS}]\label{3.4t} Let $\H$ be a right quaternionic Hilbert space and  $\{u_i\}_{i\in \NN}$ be a sequence in $V_R(\HH)$. Then, $\{u_i\}_{i\in \NN}$ is a Bessel sequence for $V_R(\HH)$ with bound $B$ if and only if the right linear operator $T: \E \to \H$ defined by
\begin{eqnarray*}
T\left(\{q_i\}_{i\in\NN}\right)= \sum_{i\in \NN}u_iq_i,\ \{q_i\}_{i\in \NN}\in\ell_2(\HH)
\end{eqnarray*}
is a well defined and  bounded operator with $\|T\|\leq\sqrt{B}$.
\end{theorem}

Let $\H$ be a right quaternionic Hilbert space and $\{u_i\}_{i\in \NN}$ be a frame for $V_R(\HH)$. Then, the \textit{(right) frame operator} $S:V_R(\HH)\rightarrow V_R(\HH)$ for the frame $\{u_i\}_{i\in \NN}$ is a right linear operator given by
\begin{eqnarray*}
S(u)&=&TT^*(u)\\
&=&T(\{\langle u_i|u\rangle\}_{i\in \NN})\\
&=&\sum_{i\in\NN}u_i\langle u_i|u\rangle,\ u\in \H.
\end{eqnarray*}

\begin{theorem}[\cite{SS}]
	Let $\H$ be a right quaternionic Hilbert space and  $\{u_i\}_{i\in \NN}$ be a frame for $V_R(\HH)$ with lower and upper frame bounds $A$ and $B$, respectively and  frame operator $S$. Then $S$ is positive, bounded, invertible and self adjoint right linear operator on $V_R(\HH)$. 
\end{theorem}

\begin{theorem}[\cite{SS}]\label{3.7t}
	Let $\H$ be a right quaternionic Hilbert space and  $\{u_i\}_{i\in \NN}$ be a frame for $V_R(\HH)$ with lower and upper frame bounds $A$ and $B$, respectively and  frame operator $S$.  Then $\{S^{-1}u_i\}_{i\in\NN}$ is also a frame for $\H$ with bounds $B^{-1}$ and ${A^{-1}}$ and right frame operator $S^{-1}$. 
\end{theorem}

\begin{theorem}[\cite{SS}]	Let $\H$ be a right quaternionic Hilbert space and $\{u_i\}_{i\in \NN}$ be a sequence in  $\H$. Then $\{u_i\}_{i\in\NN}$ is a frame for $\H$ if and only if the right linear operator $T:\E\to \H$ 
	\begin{eqnarray*}
	T(\{q_i\}_{i\in \NN})=\sum_{i\in \NN}u_iq_i, \  \ \cc{q_i}_{i\in \NN} \in \H
	\end{eqnarray*}
	is a well-defined and bounded mapping from $\ell_2(\HH)$ onto $\H$.
\end{theorem}

\section{Duals of a frame in quaternionic Hilbert spaces }
\setcounter{equation}{0}

\mathindent4em

In view of Theorem \ref{3.7t}, we have the following definition
\begin{definition}
		Let $\H$ be a right quaternionic Hilbert space and  
$\xmn$ be a  frame for $\H$ with  frame operator
$\S$. Then
\begin{itemize}[leftmargin=.3in]
\item $\cc{\S^{-1}(u_i)}_{i\in \NN}$ is called the \emph{canonical dual
frame} for the frame $\xmn$ in $\H$.
\item a sequence $\ymn \subset \H$  is called an \emph{alternate dual} for the frame  $\xmn$ in $\H$ if it satisfies
\begin{eqnarray*}
u =  \mnsum v_i\la u_i|u\ra , \ \ u \in \H.\end{eqnarray*}
\end{itemize}
\end{definition}
In view of above definition, one may observe that   canonical dual frame for a right quaternionic Hilbert space $\H$ is
also an   alternate dual frame for $\H$ and  an alternate dual of a
  frame may not be unique. In order to show their existence
we give a following example
\begin{example}
	Let $N=\cc{z_i}_{i\in\NN}$ be a Hilbert basis for a right quaternionic Hilbert space  $\H$. Then, for each $z_i, z_k\in N$, $i, k\in \NN$, we have
\begin{eqnarray*}
\langle z_i|z_k\rangle=
\begin{cases}
0,\ \text{for}\ i\neq k\\
1,\ \text{for}\ i=k.
\end{cases}
\end{eqnarray*}
Define a sequence $\cc{u_i}_{i\in \NN} \subset \nh$ by
\begin{eqnarray*}
u_{i} = z_j , \ \ i \in \cc{2j-1, 2j}, \ \ j =1,2,\cdots .
\end{eqnarray*}
Then $\xmn$ is   frame for  $\nh$. Moreover,  the  canonical dual
$\cc{\widetilde{u}_{i}}_{i\in \NN}$ of $\xmn$ is given by
\begin{eqnarray*}
\widetilde{u}_{i} = \dfrac{z_j}{2} , \ \ i \in \cc{2j-1, 2j}, \ \ j =1,2,\cdots .
\end{eqnarray*}
Next, define sequences  $\ymn$ and $\cc{w_{i}}_{i\in \NN}$ in $\nh$ by
\begin{eqnarray*}
v_{2i-1} = z_i \ \ \text{and}\ \  v_{2i} = 0  , \ \ \ i =1,2,\cdots .
\end{eqnarray*}
and
\begin{eqnarray*}
w_{2i-1} = 0 \ \ \text{and}\ \  w_{2i} = z_i  , \ \ \ i =1,2,\cdots .
\end{eqnarray*}
Then $\ymn$ and  $\cc{w_i}_{i\in \NN}$ are alternate duals for the frame $\xmn$ in $\nh$.\hfill $\Box$
\end{example}
In the following result we show that if a Bessel sequence is an   alternate dual for a given    frame in quaternionic Hilbert space $\nh$, then it becomes a frame for $\nh$ and the given   frame becomes  its   alternate dual.

\begin{theorem}
Let $\nh$ be a right quaternionic Hilbert space and $\xmn$ be a  frame for $\nh$ with lower and upper frame bounds  $A$ and $B$, respectively and $\ymn$ be a  Bessel sequence for $\nh$.
If $\ymn$ is a alternate dual  for $\xmn$ in $\nh$, then $\ymn$ is also
a frame for $\nh$. Further, $\xmn$ is a alternate dual for $\ymn$ in $\nh$.
\end{theorem}

\proof As $\xmn$ is a frame for $\nh$, therefore $\cc{\la u_i |u \ra}_{i\in \NN}\in \elii, \ \ u \in \nh.$
Also, $\ymn$ is a Bessel sequence, so by Theorem 3.2, $   \mnsum v_{i} \la  u_{i} | u\ra \ \ \text{exists}.$ So, we have
\begin{eqnarray*}
\la u|v \ra &=&  \left\la u \bigg|\mnsum v_i\la u_i|v\ra\right\ra\\
&=&  \mnsum \la u|v_{i}\ra \la u_{i}| v\ra\\
&=&\left\la \mnsum u_i \overline{\la u|v_i\ra} \bigg|v\right\ra\\
 &=&\left\la \mnsum u_i \la v_i|u\ra \bigg|v\right\ra, \ \ v\in \nh.
\end{eqnarray*}
This gives $u = \mnsum u_i\la v_i| u\ra, \  u \in \nh.$ Then we have
\begin{eqnarray*}
\pp{u}^2 &=& \sup\limits_{\pp{v}=1} |\la u
|v \ra|^2\\
&=& \sup\limits_{\pp{v}=1}   \left| \left\la \mnsum u_i\la v_i  | u\ra\bigg| v \right\ra\right|^2\\
&\le& \sup\limits_{\pp{v}=1}   \left|\mnsum \la u | v_i\ra \la  u_i | v  \ra\right|^2\\
&\le& \sup\limits_{\pp{v}=1} \bigg(  \mnsum |\la u |v_i\ra|^2\bigg)  \ \bigg( B\pp{v}^2\bigg) \\
&=& B  \mnsum |\la v_i|u\ra|^2, \ \ u\in \nh.
\end{eqnarray*}
Hence, $\ymn$ is   frame for $\nh$ with   alternate dual frame $\xmn$ in $\nh$.\endproof

In the next result, we show that among all the representations of an element\break$u \in \nh$ in terms of a  frame $\xmn$ for $\nh$ with coefficient sequence in $\elii$, the sequence $\cc{\la \S^{-1} (u_{i})|u\ra}_{i\in \NN} \in\elii$,  has the minimum $\elii$ - norm. Indeed we have the following:

\begin{theorem}
Let $\nh$ be a right quaternionic Hilbert space and $\xmn$ be a frame for $\nh$ with the  frame operator $\S$.
Fix $\tilde{u}\in \nh$, if $\tilde{u} = \mnsum u_{i}q_{i} $, for some quaternion sequence $\cc{q_{i}}_{i\in \NN} \in \elii$ then
\begin{eqnarray*}
  \mnsum|q_{i}|^2 =   \mnsum|\la \S^{-1} u_{i}| \tilde{u}\ra|^2 +
  \mnsum|\la  \S^{-1} u_{i}|\tilde{u}\ra - q_{i}|^2.
\end{eqnarray*}
In particular,   $\cc{\la  \S^{-1} u_{i} |\tilde{u} \ra}_{i \in \NN}$ has the minimal $\elii$-norm among all quaternion  sequences
$\cc{q_{i}}_{i\in \NN} \in \elii$.
\end{theorem}

\proof For each $u\in \nh$, we have
\begin{align*}
\la \S^{-1} u|u\ra &=  \left\la S^{-1}u \bigg| \mnsum u_i\la  \S^{-1} u_{i}|u\ra 
 \right\ra\\
&=   \mnsum \la \S^{-1} u|  u_{i} \ra \la \S^{-1} u_{i}|  u \ra \\
&= \left\la \cc{\la \S^{-1} u_{i}|u\ra}_{i\in \NN}| \cc{\la \S^{-1} u_{i}|u\ra}_{i\in \NN}\right\ra_{\elii}.
\end{align*}
Also
\begin{align*}
\la \S^{-1} \tilde{u}|\tilde{u} \ra &= \left\la \S^{-1} \tilde{u}\bigg|\mnsum  u_{i}q_i \right\ra\\
&= \mnsum  \la \tilde{u}|\S^{-1} u_{i}  \ra q_i\\
&= \left\la \cc{\la    \S^{-1} u_{i}|\tilde{u}\ra}_{i\in \NN}|\cc{q_{i}}_{i\in \NN}\right\ra_{\elii}.
\end{align*}
So, $\cc{q_{i}-\la  \S^{-1} u_{i}|\tilde{u}\ra}_{i\in \NN}$ is orthogonal to $\cc{\la  \S^{-1} u_{i}|\tilde{u}\ra}_{i\in \NN}$ in $\elii$.  Therefore, we have
\begin{align*}
\pp{\cc{q_i}_{i\in \NN}}_{\elii}^2 &=\pp{\cc{q_{i}-\la  \S^{-1} u_{i}|\tilde{u}\ra} +
\cc{\la  \S^{-1} u_{i}|\tilde{u}\ra}}_{\elii}^2\\
 &= \pp{\cc{q_{i}-\la  \S^{-1} u_{i}|\tilde{u}\ra}}_{\elii}^2 + \pp{\cc{\la  \S^{-1} u_{i}|\tilde{u}\ra}}_{\elii}^2.
 \end{align*}
 Thus, we have
 \begin{eqnarray*}
 \mnsum|q_{i}|^2 &=& \mnsum|\la \S^{-1} u_{i}|\tilde{u}\ra|^2 +
 \mnsum|\la\S^{-1} u_{i} |\tilde{u}\ra - q_{i}|^2. \quad\qquad\qquad\qquad\qquad\quad\qquad\Box
\end{eqnarray*}

\def\T{T}
\def\U{U}
\def\II{\mathcal{I}}
\def\cmn{\cc{q_{i}}_{i\in \NN}}

Next, we give equivalent conditions for two  frames in a quaternionic Hilbert space,  where one becomes alternate dual of the other and vice versa, in terms of their corresponding  analysis and synthesis operators. More precisely we have :

\begin{theorem}\label{t35}
Let $\nh$ be a right quaternionic Hilbert space. Let 
$\xmn$ and $\ymn$ be the frames for $\nh$ with synthesis operators $\T$
and $\U$, respectively. Then the following statements are equivalent
\begin{enumerate}[label = {\emph{(\alph*).}}]\itemsep=0em\topsep=-1em
\item $\xmn$ is an alternate  dual for $\ymn$ in $\nh$.
\item $\T\U^* =\II$.
\item $\U\T^* = \II$.
\item $\left( \T^*\U \right)^2 = \T^* \U$.
\end{enumerate}
\end{theorem}

\def\y{  \mnsum \ v_{i}q_{i}}
\proof $(a)\Rightarrow(b)$ For each $u\in\nh$, we have
\begin{eqnarray*}
u &=&   \mnsum v_{i}\la  u_{i} | u\ra \\
&=& \U \T^* (u).
\end{eqnarray*}

\noindent
$(b)\Rightarrow(c)$ Straight forward.

\noindent
$(c)\Rightarrow(d)$ We have,
\begin{eqnarray*}
\bigg( \T^*\U \bigg)^2 &=&  \T^*( \U  \T^*) \U \\
&=& \T^*\II \U \\
&=& \T^* \U.
\end{eqnarray*}

\noindent
$(d)\Rightarrow(a)$ For each $\cmn \in \elii$, we have
\begin{eqnarray}\label{e31}
\T^* \U \big(\cmn\big) = \left\{{\left\la u_i \bigg| \mnsum \ v_{i} q_i\right\ra}\right\}_{i\in \NN}.
\end{eqnarray}
This gives
\begin{align*}
\U \T^* \U \left(\cmn\right) = \mnsum v_i  {\left\la u_i \bigg| \mnsum \ v_{i} q_i\right\ra}.
\end{align*}
Therefore, we have
\begin{align}\label{e32}
(\T^*  \U)^2 \big(\cmn\big)  = \left\{\left\la u_i \bigg| \mnsum v_i{\left\la u_i \bigg| \mnsum \ v_{i} q_i\right\ra}\right\ra\right\}_{i\in \NN}.
\end{align}
So (\ref{e31}) and (\ref{e32}) gives
\begin{eqnarray*}
\y = \mnsum v_i\left\la u_i \bigg| \mnsum \ v_{i} q_i\right\ra.
\end{eqnarray*}
Again since $\ymn$ is an  frame for $\nh$, therefore $ \U$ is onto. Thus we have
\begin{eqnarray*}
 v = \mnsum v_i\la  u_{i}|v\ra, \ \ v\in \nh.
\end{eqnarray*}
Hence $\ymn$ is an alternate dual frame for $\xmn$ in $\nh$. \endproof

\def\Tm{T^*}
\def\Ti{\widetilde{T}}
\def\Tmi{\widetilde{T}^*}

Next, we give a result concerning a relationship between the analysis
operator and  the canonical dual
 of a frame  in a right quaternionic Hilbert space. Further, an expression for the pseudo inverse of the  synthesis operator of a frame in terms of its  canonical dual  is given.

\begin{theorem}
Let $\nh$ be a right quaternionic Hilbert space and 
$\xmn$ be a  frame for $\nh$ with  the  frame operator $\S$. Let
$\T$ and $\Tm$ be the  synthesis operator  and  the analysis  operator, respectively  for $\xmn$ and
$\Ti$ and $\Tmi$ be the  synthesis operator  and the analysis  operator, respectively  for the canonical dual frame $\cc{\S^{-1} (u_i)}_{i\in\NN}$. Then,
\begin{enumerate}[label = (\alph*),leftmargin=.3in]\itemsep0em \topsep0em
\item range $(\Tm)$ = range $(\Tmi)$.
\item the pseudo inverse $\big( \T \big)^{\dag}$ of the  synthesis operator
$\T$ is $\Tmi$, i.e.
\begin{eqnarray*}
\big( \T \big)^{\dag} u = \big\{\la \S^{-1} (u_{i})|u\ra\big\}_{i \in \NN}, \ \ u\in \nh.
\end{eqnarray*}
\end{enumerate}
\end{theorem}

\proof $(a)$ For each $u\in\nh$, we have
\begin{eqnarray*}
\big( \Tmi \big) u &=& \big\{\la \S^{-1} (u_{i} )|u\ra\big\}_{i\in\NN}\\
&=& \big\{\la\ u_i|S^{-1}u\ra\big\}_{i\in\NN}\\
&=& \Tm \S^{-1} u.
\end{eqnarray*}
Since $\S$ is a topological isomorphism, it follows that $range \ (\Tm) = range \ (\Tmi)$.
\medskip

\noindent
$(b)$ As $\ker (\T)^\perp = \ range (\Tm) = \ range (\Tmi)$. Therefore,
$ \T|_{{\ker(\T)}^{\perp}} : range(\Tm)\to \nh$
is a topological isomorphism. Therefore the pseudo inverse $\big( \T \big)^{\dag}$ is
$\bigg( \T|_{{\ker(\T)}^{\perp}}\bigg)^{-1}$. Further,
\begin{eqnarray*}
\bigg( \T|_{{\ker(\T)}^{\perp}}\bigg) \Tmi u &=&
\T \Tmi u \\
&=&  \mnsum u_i\la \S^{-1} u_{i}|u \ra \\
&=& u ,\ \ u\in \nh.
 \end{eqnarray*}
Furthermore, for each $q\in range \ (\Tm) = range \ (\Tmi)$, there exist $u\in \nh$
such that $q = \Tmi u $. So, by Theorem \ref{t35}, we  have
\begin{eqnarray*}
\Tmi \bigg( \T|_{{\ker(\T)}^{\perp}}\bigg) q &=& \Tmi\T\Tmi u \\
&=& \Tmi u.
\end{eqnarray*}
Hence $\Tmi = \bigg( \T|_{{\ker(\T)}^{\perp}}\bigg)^{-1}.$ \endproof

\def\psxmn{\cc{\mathcal{P}\S^{-1} (u_{i})}_{i\in \NN}}
\def\pxmn{\cc{\mathcal{P}u_{ i}}_{i\in \NN}}

Finally in this section, with a  given  frame and its  canonical dual
frame for a quaternionic Hilbert space, we characterize   frame and its
  canonical dual frame for a given closed subspace of the quaternionic Hilbert
space.
\begin{theorem}
Let $\nh$ be a quaternionic Hilbert space and
$\xmn$ be a   frame for $\nh$ with  frame operators $\S$.
Let $\mathcal{P}$ be the orthogonal projection of $\nh$ onto a closed subspace $\mathcal{M}$ of $\nh$. Then
\begin{enumerate}[label = (\alph*),leftmargin=.3in] \itemsep0em \topsep0em
\item $\pxmn$ is  frame for $\mathcal{M}$ with the
same frame bounds as $\xmn$ and $\psxmn$ is an  alternative
dual frame for $\pxmn$.
\item $\psxmn$ is  canonical dual frame for $\pxmn$ if and only if
$\mathcal{P} \S = \S \mathcal{P}$.
\end{enumerate}
\end{theorem}

\proof $(a)$ For each $u\in\mathcal{M}$, we have
\begin{eqnarray*}
 \mnsum |\la \mathcal{P} u_{ i}|u \ra|^2 =   \mnsum |\la u_{i}| u  \ra|^2
\end{eqnarray*}
Therefore, $\pxmn$ is   frame for $\mathcal{M}$ with same frame bounds
as that of $\xmn$ and by the similar argument, $\psxmn$ is also  a  frame for $\mathcal{M}$ with
frame bounds inverse of $\xmn$. Further, for each $u\in \mathcal{M}$
\begin{eqnarray*}
 \mnsum \mathcal{P}u_i\la \mathcal{P}\S^{-1} u_{i}|u \ra  & =&
\mathcal{P}\bigg(\mnsum u_i\la \S^{-1} u_{i}|\mathcal{P} u\ra\bigg)\\
&=& u.
\end{eqnarray*}
Hence, the result follows.
\medskip

\noindent
$(b)$ Let $V$ be the  frame operator for $\pxmn$ as a frame for $\mathcal{M}$.
Then, we have
\begin{eqnarray*}
V ^{-1} \mathcal{P}(u_{i}) = \mathcal{P} \S^{-1}(u_{i}), \ \ i \in \NN.
\end{eqnarray*}
This gives
\begin{eqnarray*}
V ^{-1} \mathcal{P}(u) = \mathcal{P} \S^{-1}(u), \ \ \ u\in \nh.
\end{eqnarray*}
Since $\S$ and $V$ are topological isomorphisms on $\nh$ and $\mathcal{M}$, respectively, therefore
we have
\begin{eqnarray*}
\mathcal{P}\S &=& VV^{-1} \mathcal{P}\S \\
&=& V \mathcal{P}\S^{-1} \S \\
&=& V \mathcal{P}.
\end{eqnarray*}
Since $\mathcal{P}, \ \S $ and $V$ are self-adjoint, so we have $\S\mathcal{P} = \mathcal{P} V$.
Hence, we have $\S\mathcal{P} = \mathcal{P}\S$.
\medskip

\noindent
Conversely,  we have
\begin{eqnarray*}
V  \mathcal{P}u&=&   \mnsum \mathcal{P}u_{i}\la  \mathcal{P}u_{i}|\mathcal{P}u\ra  \\
&=&   \mathcal{P} \bigg( \mnsum u_{i}\la u_{i}|\mathcal{P}^2u \ra \bigg)\\
&=& \mathcal{P} \bigg(   \mnsum u_{i}\la u_{i}|\mathcal{P}u \ra \bigg)\\
&=& \mathcal{P}\S\mathcal{P}u\\
&=&\S\mathcal{P}u, \ \ \ u\in \nh.
\end{eqnarray*}
Therefore, for each $u\in \mathcal{M}$, we have
\begin{eqnarray*}
Vu &=& V \mathcal{P}u\\
&=& \S \mathcal{P}u\\ &=& \S u.
\end{eqnarray*}
This gives $V  = \S|_{\mathcal{M}}$. Thus $\S$ maps $\mathcal{M}$ bijectively onto
itself. Therefore we have
\begin{eqnarray*}
V^{-1}\mathcal{P}(u_{i}) &=& \big(\S|_{\mathcal{M}}\big)^{-1} \mathcal{P} u_{i} \\
&=&\S^{-1} \mathcal{P}(u_{i}) \\
&=&\mathcal{P}\S^{-1} (u_{i}).
\end{eqnarray*}
Hence the   canonical dual frame for $\pxmn$ is $\psxmn$ in $\nh$. \endproof

\section{Application}
In this section, as an application of the  canonical dual of a given frame, we give the orthogonal projection of $\elii$ onto the range of the analysis operator of the  frame, in terms of the elements of the  frame and  its canonical dual.  Apart from this, we also give an expression for the orthogonal projection in terms of operators related to the   frame and its   canonical dual  in a right quaternionic Hilbert space.
\begin{theorem}
Let $\nh$ be a right quaternionic Hilbert space and $\xmn$ be a frame for $\nh$ with the synthesis operator $T$ and the frame operator $S$. Let $\tilde{T}$ and $\tilde{S}$ be the synthesis operator and the frame operator for the  canonical dual frame $\lbrace S^{-1} u_{i}\rbrace_{i\in\NN}$ of $\cc{u_i}_{i\in \NN}$. Then the orthogonal projection $Q$ of $\elii$ onto the range of $T^*$ is given by
\begin{eqnarray*}
Q\bigg(\lbrace q_{i}\rbrace_{i\in\NN}\bigg)=\bigg\lbrace \bigg\la u_{j}\bigg|\sum\limits_{i=1}^{\infty}S^{-1}u_{i}q_{i} \bigg\ra\bigg\rbrace_{j \in\NN}.
\end{eqnarray*}
 Also,
\begin{eqnarray*}
Q = T^*\tilde{T} 
=\T^*\S^{-1}\T.
\end{eqnarray*}
\end{theorem}
\proof
For $q=\lbrace q_{i}\rbrace_{i\in\NN} \in \elii$,
\begin{eqnarray*}
u&=&S^{-1}\S u\\
&=&\sum\limits_{i=1}^{\infty}\S^{-1}u_{i}\la u_{i}|u \ra , \ \ u\in \nh.
\end{eqnarray*}
Moreover, we have $  \T(q)= \mnsum u_{i}q_{ i}\ \  \text{and} \ 
\T^*u=\lbrace \la u_{i}|u\ra\rbrace_{i\in \NN}.$ 
It is sufficient to show that $Q$ is the identity on $range\ {T^*}$ and  is zero on $(range\ {T^*})^\perp$=$\ker{\T}$. By definition
\begin{eqnarray*}
Q\big(T^*u\big)&=&Q\left(\{ \la u_i|u \ra\}_{i\in\NN}\right)\\
&=&\bigg\lbrace\bigg\la u_{j} \bigg|\sum\limits_{i=1}^{\infty}S^{-1}u_{i}\la u_{i}|u\ra \  \bigg\ra\bigg\rbrace_{j\in\NN}
\\
&=&\lbrace\la u_{ j}|u \ra \rbrace_{j\in\NN}\\
&=&T^*u, \ \ u \in \nh.
\end{eqnarray*}
Again, for $q=\lbrace q_{i}\rbrace_{j\in\NN}\in (range \ T^*)^\perp=\ker \ \T$
\begin{eqnarray*}
Q(q)&=& Q\bigg(\lbrace q_{i}\rbrace_{i\in\NN}\bigg)\\
&=& \bigg\lbrace \bigg\la u_{j}\bigg|\sum\limits_{i=1}^{\infty}S^{-1}u_{i}q_{i} \bigg\ra\bigg\rbrace_{j\in\NN}\\
&=&\bigg\lbrace\bigg\la u_{j}\bigg|S^{-1} \sum\limits_{i=1}^{\infty}u_{i}q_{i} \bigg\ra \bigg\rbrace_{j\in\NN}\\
&=&\left\lbrace\big\la u_j |S^{-1}T(q)\big\ra \right\rbrace_{j\in\NN}\\ &=&0.
\end{eqnarray*}
 Thus, Q  is the orthogonal projection of
 $\elii$ onto the range of $T^*$.
Further, we have
\begin{eqnarray*}
T^*S^{-1}T(q)&=&T^*\tilde{T}(q)\\
&=&T^*\bigg(\sum\limits_{i=1}^{\infty}S^{-1}u_{i}q_{i}\bigg)\\
&=&\bigg\lbrace\bigg\la u_j \bigg|\sum\limits_{i=1}^{\infty}S^{-1}
u_{i}q_{i}\bigg\ra\bigg\rbrace_{j\in\NN}\\
&=&Q(q).
\end{eqnarray*}
This completes the proof. \endproof

\end{document}